\documentclass[11pt]{article}
\usepackage{latexsym,amssymb,amsthm,amsmath,multirow}
\usepackage{amssymb}

\def\<{\left < }
\def\>{\right >}
\def\({\left ( }
\def\){\right )}
\def\e{\eqref}
\def\E4{\mathbb E^4 }
\def\g{\gamma }

\def\g{\gamma}

\begin{document}

\vbox{\hbox{\small Proceedings of the Conference RIGA 2011}
\hbox{\small {\it Riemannian Geometry and Applications}}
\hbox{\small Bucharest, Romania}}
\vskip 1.5truecm

\setlength{\textheight}{19cm}
\setlength{\textwidth}{12.5cm}

\centerline{\large{\bf ON WINTGEN IDEAL SURFACES}}\medskip
\centerline{\bf Bang-Yen Chen}

\bigskip

\begin{abstract}  Wintgen proved in \cite{W} that  the Gauss curvature $K$ and the normal curvature $K^D$ of a surface  in the Euclidean 4-space $\E4$ satisfy $$K+|K^D|\leq H^2,$$ where $H^2$ is the squared mean curvature. A surface $M$ in $\E4$ is called a {Wintgen ideal} surface if it satisfies the equality case of the inequality identically. Wintgen ideal surfaces in $\E4$ form an important family of surfaces; namely, surfaces with circular ellipse of curvature.
In this paper, we provide a brief survey on some old and recent results on Wintgen ideal surfaces and more generally Wintgen ideal submanifolds in definite and indefinite real space forms. 

\vskip.04in
\noindent 2000  {\it Mathematics Subject Classification}: {Primary: 53A05; Secondary  53C40, 53C42}
\vskip.04in

\noindent {\it Keywords}: {Gauss curvature, normal curvature, squared mean curvature,  Wintgen ideal surface, superminimal surface, Whitney sphere.}

\end{abstract}

\section{Introduction}

For surfaces $M$ in a Euclidean 3-space $\mathbb E^3$, the {\it Euler inequality} 
\begin{align}\label{1.0}K\le H^2,\end{align}
 whereby $K$ is the intrinsic {\it Gauss curvature} of $M$ and $H^2$ is the extrinsic {\it squared mean curvature} of $M$ in $\mathbb E^3$, at once follows from the fact that $$K=k_1k_2,\;\; H=\frac {1}{2}(k_1+k_2),$$ whereby $k_1$ and $k_2$ denote the {\it principal curvatures} of $M$ in $\mathbb E^3$. And, obviously, $K=H^2$ everywhere on $M$ if and only if the surface $M$ is {\it totally umbilical} in $\mathbb E^3$, i.e. $k_1=k_2$ at all points of $M$, or still, by a {\it theorem of Meusnier}, if and only if $M$ is a part of a {\it plane} $E^2$ or of a {\it round sphere} $S^2$ in $\mathbb E^3$. 

Consider an isometric immersion $\psi:M\to \tilde M^4$ of a surface $M$ into a Riemannian 4-manifold $\tilde M^4$,   the {\it ellipse of curvature} at a point $p$ of $M$ is defined as \begin{align}\label{1.2}\mathcal E_p = \{h(X,X)\,|\, X\in T_pM, \|X\|=1\},\end{align} where $h$ is the second fundamental form of $M$ in $\tilde M^4$.

In  1979, P. Wintgen \cite{W} proved a basic relationship between the intrinsic {\it Gauss curvature} $K$, the extrinsic {\it normal curvature} $K^D$, and {\it squared mean curvature} $H^2$ of any surface $M$ in a Euclidean 4-space $\mathbb E^4$; namely,  \begin{align} \label{1.1}K+|K^D|\leq H^2,\end{align} 
with the equality holding if and only if the {curvature ellipse}  is a {circle}. Wintgen's inequality was generalized to surfaces in 4-dimensional real space forms in \cite{GR}.
  A similar inequality holds for surfaces in pseudo-Euclidean 4-space $\mathbb E^{4}_{2}$ with neutral metric \cite{c2010,c6}.

Following L. Verstraelen et al. \cite{dpv,PV}, we call a surface $M$ in $\E4$ {\it Wintgen ideal} if it satisfies the equality case of Wintgen's inequality identically. Obviously, Wintgen ideal surfaces in $\mathbb E^4$ are exactly  superminimal surfaces.

In this article,  we provide a brief survey on some old and some recent results on Wintgen ideal surfaces; and more generally, Wintgen ideal submanifolds in definite and indefinite real space forms. Some related results are also presented in this paper.

\pagestyle{myheadings}

\markboth{B.-Y. Chen}{WINTGEN IDEAL SURFACES}

 \section{Some known results on superminimal surfaces}
    
\subsection{$R$-surfaces}

A surface $\psi:M\to \tilde M^4$ is  {\it superminimal} if and only if, at each point $p\in M$, the ellipse of curvature $\mathcal E_p$ is a circle with center at the origin $o$ (see \cite{F1}). Simple examples of superminimal surfaces in  the Euclidean 4-space $\mathbb E^4$ are $R$-surfaces, i.e., graphs of holomorphic functions: 
\begin{align}\{(z,f(z)):z\in U\},\end{align}
  where $U\subset {\mathbb C} \approx {\bf R}^2$ is an open subset of the complex plane and $f$ is a holomorphic function.

When the ambient space $\tilde M^4$ is a space of constant curvature, O. Bor\r{u}vka \cite{B} proved in 1928 that the family of superminimal immersions $\psi:M\to \tilde M^4$ depends (locally) on two holomorphic functions.

 \subsection{Isoclinic surfaces}
For an oriented plane $E$ in $\mathbb E^4$, let $E^\perp$ denote the orthogonal complement with the orientation given by the condition $$E\oplus E^\perp=\mathbb E^4.$$ Two oriented planes $E,F$ are called {\it oriented-isoclinic} if either 
\vskip.05in

\noindent (a) $E=F^\perp$ (as oriented planes) or 
\vskip.05in

\noindent (b) the projection $pr_F:E\to F$ is a non-trivial, conformal map preserving the orientations.
 \vskip.05in
 
  Consider an oriented surface $\psi:M\to \tilde M^4$. If $\gamma$ is a curve in $\tilde M^4$, denote by $\tau_\gamma$ the parallel displacement along $\gamma$ in the tangent bundle $T\tilde M^4$.   The surface $M^2$ is called a {\it negatively oriented-isoclinic surface} if, for every curve $\gamma$ in $M$ from $x$ to $y$, the planes $\tau_{\psi\circ \gamma}(T_{\psi(x)}M)$ and  $T_{\psi(y)}M$ are negatively oriented isoclinic planes in $T_{\psi(y)}\tilde M^4$.

S. Kwietniewski proved  in his 1902 dissertation at Z\"urich \cite{Kw} that a surface in $\mathbb E^4$ is superminimal if and only if it is negatively oriented-isoclinic.  
Th.   Friedrich in 1997 extended this result  for surface in $\tilde M^4$.

\subsection{Representation}
   In 1982, R. Bryant \cite{Br} studied a superminimal immersion of a Riemann surface $M$ into $S^{4}$ by lifting it to $CP^{3}$, via the twistor map $$\pi:CP^{3}\to S^{4}$$ of Penrose. The lift is a holomorphic curve, of the same degree as that of the immersion, which is horizontal with respect to the twistorial fibration; moreover, the lift is a holomorphic curve in $CP^{3}$ satisfying the differential equation 
   \begin{align} z_{0}dz_{1}-z_{1}dz_{0}+z_{2}dz_{3}-z_{3}dz_{2}=0.\end{align}
   
Setting $$z_{0}=1, \;\; z_{1}+z_{2}z_{3}=f,\;\; z_{2}=g,$$one can solve for $z_{1},z_{2},z_{3}$ in terms of the meromorphic functions $f$ and $g$, which serves as a kind of  
Weierstrass representation. 
   Via this, R. Bryant showed the existence of a superminimal immersion from any compact Riemann surface $M$ into the 4-sphere $S^{4}$.

M. Dajczer and R. Tojeiro  established in \cite{DT} a representation formula for superminimal surfaces in $\mathbb E^4$ in terms of pairs $(g,h)$ of conjugate minimal surfaces in $\mathbb E^4$.

\subsection{Twister space}
On an oriented Riemannian 4-manifold $\tilde M^4$,  there exists an $S^{2}$-bundle $Z$, called the {\it twistor space} of $\tilde M^4$, whose fiber over any point $x\in \tilde M^4$ consists of all almost complex structures on $T_{x}\tilde M^4$ that are compatible with the metric and the orientation. It is known that there exists a one-parameter family of metrics $g^{t}$  on $Z$, making the projection 
\begin{align} Z\rightarrow \tilde M^4\end{align} into a Riemannian submersion with totally geodesic fibers.
       
Th.  Friedrich proved in 1984 that  superminimal surfaces are characterized by the property that the lift into the twistor space is holomorphic and horizontal.

\subsection{Central sphere congruence}

The  central sphere congruence of a surface in Euclidean space is the family of 2-dimensional spheres that are tangent to the surface and have the same mean curvature vector as the surface at the point of tangency.

 In 1991, B. Rouxel \cite{R} proved the following results:
\vskip.1in
   
\noindent {\bf Theorem 2.1.} {\it If the ellipse of curvature of a surface in $\mathbb E^{4}$ is a circle, then the surface of centers of the harmonic spheres is a minimal surface of $\mathbb E^{4}$.}
\vskip.04in

\vskip.1in
   
\noindent {\bf Theorem 2.2.} {\it  If $M$ is a surface of $\mathbb E^{4}$ with circular ellipse of curvature and if the harmonic spheres of $M$ have a common fixed point, then $M$ is a conformal transform of a superminimal surface of $\mathbb E^{4}$.}

\vskip.1in
   
\noindent {\bf Theorem 2.3.} {\it The surface of centers of such sphere congruence is a minimal surface.}

\subsection{Ramification divisor} 

   Let $M$ be a compact Riemann surface of genus $g$ and let $\phi:M\to CP^{1}$ be a holomorphic map of degree $d$. A point $x\in M$ is a {\it ramification point} of $\phi$ if $d\phi(x)=0$, and its image $\phi(x)\in CP^{1}$ is called a {\it branch point} of $\phi$.
   
  By the Riemann-Hurwitz Theorem the number of branch points of $\phi$ (counting multiplicities) is $2g+2d-2$.
   
The {\it ramification divisor} of $\phi$ is the formal sum 
\begin{align}\sum_{i}a_{i}p_{i},\end{align}
 where $p_{i}$ is a ramification point of $\phi$ with multiplicity $a_{i}$, and where the sum is taken over all ramification points of $\phi$. 
  Let $Ram(\phi)$  denote the ramification divisor of $\phi$.

If we put  \begin{align}f_{1}=\frac{z_{1}}{z_{0}},\;\; f_{2}=\frac{z_{3}}{z_{2}},\end{align} then $f_{1}$ and $f_{2}$ are known of degree $d$  satisfying $ram(f_{1})=ram(f_{2})$, where $ram(f)$ is the {\it ramification divisor}  of the meromorphic function $f$. 
   
This provides a method for constructing the {\it moduli space} ${\mathcal M}_{d}(M)$ of horizontal holomorphic curves of degree $d$ for a Riemann surface  $M$ in  $S^{4}$.

 For $M=S^{2}$, B. Loo  proved \cite{Loo}  that the moduli space ${\mathcal M}_{d}(M)$ is connected and it has dimension $2d+4$.

\subsection{Riemann surfaces of higher genera}

 By applying algebraic geometry,  Chi and  Mo studied in \cite{CM} the moduli space over superminimal surfaces of higher genera. In particular, they proved the following 5 results:

\vskip.1in
   
\noindent {\bf Theorem 2.4.} {\it  Let  $M$ be a Riemann surface of genus $g\geq 1$.	Then all the branched superminimal immersions of degree $d<5$ from $M$ into $S^{4}$ are totally geodesic.}
\vskip.1in
   
\noindent {\bf Theorem 2.5.} {\it  Let  $M$ be a Riemann surface of genus $g\geq 1$.	Then $M$ admits a non-totally geodesic branched superminimal immersions of degree $6$  into $S^{4}$ if and only if $M$ is a hyper-elliptic surface, i.e., it is an elliptic fibration over an elliptic curve.}

\vskip.1in
   
\noindent {\bf Theorem 2.6.} {\it   Let $M$ be a hyper-elliptic surface of genus $g>3$. Then non-totally geodesic branched superminimal immersions of degree 6 from $M$ into $S^{4}$ are the pullback of non-totally geodesic branched superminimal spheres of degree 3 via the branched double covering of $M$ onto $CP^{1}$.} 
\vskip.1in
   
\noindent {\bf Theorem 2.7.} {\it  Let $M$ be a Riemann surface of genus $g\geq 2$ ($g=1$, respectively).	If  $d>5g+4$,  ($d\geq 6$, respectively), then there is a non-totally geodesic branched superminimal immersion of degree $d$ from $M$ into $S^{4}$. The immersion is generically one-to-one.}
 \vskip.1in
   
\noindent {\bf Theorem 2.8.} {\it  Let $M$ be a  Riemann surface of genus $g\geq 1$. If the degree $d$ of a superminimal immersion of $M$ in $S^{4}$ satisfies $d\geq 2g-1$, then the dimension of  the moduli space ${\mathcal M}_{d}(M)$ is between $2d-4g +4$ and $2d-g +4$, where the upper bound is achieved by the totally geodesic component. }

\section{Wintgen's inequality}

We recall the following result of P. Wintgen \cite{W}. 

\vskip.1in
\noindent {\bf Theorem 3.1. }\label{T:3.1} {\it Let $M$ be a   surface in Euclidean 4-space $\mathbb E^4$. Then we have
 \begin{align} \label{3.1}& H^2\geq K+|K^D|\end{align}
 at every point in $M$. Moreover, we have
 \vskip.05in
 
\noindent {\rm (i)}  If $K^D\geq 0$ holds  at a point $p\in M$, then the equality sign of \e{3.1} holds at $p$ if and only if, with respect to some suitable orthonormal frame $\{e_1,e_2,e_3,e_4\}$ at $p$,  the shape operator at $p$ satisfies
   \begin{align}& \label{3.2} A_{e_3}=\begin{pmatrix} \mu+2\gamma &0\\0 & \mu \end{pmatrix}, \; A_{e_4}=\begin{pmatrix} 0 &\gamma\\\gamma & 0 \end{pmatrix}.\end{align}
   
 \vskip.05in
 
\noindent {\rm (ii)} If $K^D< 0$ holds  at  $p\in M$, then the equality sign of \e{3.1} holds at $p$ if and only if, with respect to some suitable orthonormal frame $\{e_1,e_2,e_3,e_4\}$ at $p$,  the shape operator at $p$ satisfies
   \begin{align}& \label{3.3} A_{e_3}=\begin{pmatrix} \mu-2\gamma &0\\0 & \mu\end{pmatrix}, \; A_{e_4}=\begin{pmatrix} 0 &\gamma\\\gamma & 0 \end{pmatrix}.\end{align} }

\section{Wintgen ideal surfaces in $\mathbb E^4$}

In this and the next sections we present some recent results on Wintgen ideal surfaces.

\vskip.1in
\noindent {\bf Proposition 4.1.}\label{P:4.1} {\it Let $M$ be a  Wintgen ideal surface in $\E4$. Then $M$ has  constant mean curvature and constant Gauss curvature if and only if  $M$ is totally umbilical.}
\vskip.1in

The following results classifies Wintgen ideal surfaces in $\mathbb E^{4}$ with equal Gauss and normal curvatures.

\vskip.1in
\noindent {\bf Theorem 4.1.}\label{T:4.1} {\it Let $\psi:M\to \E4$ be a Wintgen ideal surface in $\E4$. Then  $|K|=|K^D|$ holds identically  if and only if one of the following four cases occurs:

\vskip.04in 
\noindent {\rm (1)}  $M$ is an open portion of a totally geodesic plane in $\E4$.

\vskip.04in
\noindent {\rm (2)} $M$ is a complex curve lying fully in ${\mathbb C}^2$, where ${\mathbb C}^2$ is the Euclidean 4-space $\E4$ endowed with some orthogonal almost complex structure.

\vskip.04in 
\noindent {\rm (3)} Up to dilations and rigid motions on the Euclidean 4-space $\E4$, $M$ is an  open portion of the Whitney sphere defined by
\begin{align}\notag &\psi(u,v)=\frac{\sin u }{1+\cos^2 u}\Big(\sin v,\cos v,\cos u\sin v,\cos u\cos v
 \Big).\end{align} 

\vskip.04in 
\noindent {\rm (4)} Up to dilations and rigid motions of the Euclidean 4-space  $\E4$, $M$ is a surface with $K=K^D=\frac{1}{2}H^2$ defined by}
\begin{equation}\begin{aligned}\notag  &\hskip.1in \psi(x,y)=\text{$\frac{2 \sqrt{y}}{5}$}\sqrt{\cos x} \cos\(\text{$\frac{x}{2}$}\)\cos(\ln y) \cos \(\text{$\frac{1}{2}$}\tanh^{-1}\(\tan \text{$\frac{x}{2}$}\)\)  \\&\hskip.1in \times  \Bigg(\hskip-.02in\tan \(\text{$\frac{1}{2}$}\tanh^{-1}\(\tan \text{$\frac{x}{2}$}\)\)(2-\tan(\ln y))  +\tan\(\text{$\frac{x}{2}$}\)(1+2\tan(\ln y)),
\\&\hskip.2in  \tan \(\text{$\frac{1}{2}$}\tanh^{-1}\(\tan \text{$\frac{x}{2}$}\)\)(1+2 \tan(\ln y)) - \tan\(\text{$\frac{x}{2}$}\)(2-\tan(\ln y)),
\\&\hskip.3in  
  \tan\(\text{$\frac{x}{2}$}\)\tan \(\text{$\frac{1}{2}$}\tanh^{-1}\(\tan \text{$\frac{x}{2}$}\)\)(1+ 2\tan(\ln y))+ \tan(\ln y)-2,
\\&\hskip.4in   \tan\(\text{$\frac{x}{2}$}\)\tan \(\text{$\frac{1}{2}$}\tanh^{-1}\(\tan \text{$\frac{x}{2}$}\)\)(\tan(\ln y)-2) -2\tan(\ln y)-1\hskip-.01in \Bigg).
\end{aligned}\end{equation}
\vskip.1in

According to I. Castro \cite{Ca}, up to rigid motions and dilations of ${\mathbb C}^2$ the Whitney sphere is the only compact orientable Lagrangian superminimal surface in ${\mathbb C}^2$.
\vskip.1in

\noindent {\bf Remark 4.1.} In order to prove Theorem 4.1  we have solved the following fourth order differential equation:
 \begin{equation}\begin{aligned}\label{5.61} & p^{(4)}(x)-2(\tan x) p'''(x)+\(1+\text{$\frac{5}{8}$}\sec^2 x\) p''(x)\\&\hskip.6in
 +\(\text{$\frac{5}{8}$}\sec^2x -2\)(\tan x ) p'(x)+\text{$\frac{185 }{256}$}(\sec^4 x)p(x)=0.
 \end{aligned}\end{equation}
to obtain the following exact solutions: 
  \begin{equation}\begin{aligned}\label{5.62} & p(x)= \sqrt{\cos  x}\left\{ \!\(c_1 \cos\!  \text{$\(\frac{x}{2}\)$}+c_2\sin  \!\text{$\(\frac{x}{2}\)$}\)\cos \(\text{$\frac{1}{2}$} \tanh^{-1}\(\tan\! \text{$\(\frac{x}{2}\)$}\!\)\!\)\right. \\&\left. \hskip.4in + \(c_3 \cos \text{$\(\frac{x}{2}\)$}+c_4\sin  \text{$\(\frac{x}{2}\)$}\)\sin \(\text{$\frac{1}{2}$} \tanh^{-1}\(\tan \text{$\(\frac{x}{2}\)$}\)\)\right\}
 \end{aligned}\end{equation}

\section{Wintgen ideal surfaces in $\mathbb E^{4}_{2}$}

For space-like oriented surfaces in a 4-dimensional indefinite real space form $R^4_2(c)$ with neutral metric, one has the following Wintgen type inequality (cf. \cite{c2010,c6,CS}).

\vskip.1in
\noindent {\bf Theorem 5.1.}\label{T:5.1} {\it  Let $M$ be an oriented space-like surface in a 4-dimensional indefinite space form $R^4_2(c)$ of constant sectional curvature $c$ and with index two. Then we have
 \begin{align} \label{5.1}& K+K^D\geq \<H,H\>+c \end{align}
 at every point. 
 
 The equality sign of \e{5.1} holds at a point $p\in M$ if and only if, with respect to some suitable orthonormal frame $\{e_1,e_2,e_3,e_4\}$,  the shape operator at $p$ satisfies
   \begin{align}& \label{5.2} A_{e_3}=\begin{pmatrix} \mu+2\gamma &0\\0 & \mu \end{pmatrix}, \; A_{e_4}=\begin{pmatrix} 0 &\gamma\\\gamma & 0 \end{pmatrix}.\end{align}}
   \vskip.1in

As in surfaces in 4-dimensional real space forms, we call a surface in $R^{4}_{2}(c)$  Wintgen ideal if it satisfies the equality case of \e{5.1} identically.

\vskip.1in
\noindent {\bf Theorem 5.2.}\label{T:5.2} {\it Let $M$ be a Wintgen ideal surface in a neutral pseudo-Euclidean $4$-space $\mathbb E^4_2$. Then $M$ satisfies $|K|=|K^D|$ identically if and only if, up to dilations and rigid motions, $M$ is one of the following three types of surfaces:

\vskip.04in
\noindent {\rm (i)}   A space-like complex curve in ${\bf C}^2_1$, where ${\bf C}^2_1$ denotes $\mathbb E^4_2$ endowed with some orthogonal complex structure;

\vskip.04in
\noindent {\rm (ii)}  An open portion of a non-minimal  surface defined by
\begin{align}\notag & \text{$ \sec^2\! x\Big(\hskip-.02in \sin x \sinh y,\sqrt{2-\sin^2\! x} \cosh y ,\sin x \cosh y, \sqrt{2-\sin^2 \!x} \sinh y \Big)$};\end{align}

\vskip.04in
\noindent {\rm (iii)} An open portion of a non-minimal surface defined by
 \begin{equation}\begin{aligned}\notag  &\hskip.0in \text{$\frac{\cosh x}{6\sqrt{2y}}$}\! \text{ $\(\! 6\sqrt{2}\sqrt{\!\sqrt{2}\!+\! (1\!-\!2\tanh x)\sqrt{1\!+\!\tanh x}}+y^2\sqrt{\!\sqrt{2}\!+\!\sqrt{1\!+\!\tanh x}},\right.$}
 \\& \hskip.7in \text{ $ 6\sqrt{2}\sqrt{\!\sqrt{2}\!+\! (2\tanh x\!-\!1)\sqrt{1\!+\!\tanh x}}$} \\ &\hskip.98in + \text{$ y^2 \sqrt{\!\sqrt{2}\!+\!\sqrt{1\!+\!\tanh x}}\(\! \sqrt{2}\cosh x\sqrt{1+\tanh x}-e^x\)$},
\\&\hskip.5in \text{$ 6\sqrt{2}\sqrt{\!\sqrt{2}\!+\! (1\!-\!2\tanh x)\sqrt{1\!+\!\tanh x}}  - y^2 \sqrt{\!\sqrt{2}\!+\!\sqrt{1\!+\!\tanh x}}$},
\\&\hskip.7in \text{ $  6\sqrt{2}\sqrt{\!\sqrt{2}\!+\! (2\tanh x\!-\!1)\sqrt{1\!+\!\tanh x}}\,$}
\\&\left.  \hskip.9in - \text{$y^2 \sqrt{\!\sqrt{2}\!+\!\sqrt{1\!+\!\tanh x}}\(\sqrt{2}\cosh x\sqrt{1\!+\!\tanh x}-e^x\)$}\! \). \end{aligned}\end{equation}}
\vskip.1in

\section{Surfaces with null normal curvature in $\mathbb E^4_2$.} 

The following theorem of Chen and Suceav\u{a} from \cite{CS} classifies surfaces with null normal curvature in $\mathbb E^4_2$.

\vskip.1in
\noindent {\bf Theorem 6.1.}\label{T:6.1} {\it Let $M$ be a space-like  surface in the pseudo-Euclidean 4-space $\mathbb E^4_2$. If $M$ has constant mean and Gauss curvatures and null normal curvature, then $M$ is congruent to an open part of one of the following six types of surfaces:
\vskip.04in 

\noindent {\rm (1)} A  totally geodesic plane in $\mathbb E^4_2$ defined by $(0,0,x,y)$;
\vskip.04in 

\noindent  {\rm (2)}  a totally umbilical hyperbolic plane $H^2(-\frac{1}{a^2})\subset \mathbb E^3_1\subset \mathbb E^4_2$  given by $$\big(0, a\cosh u,a \sinh u \cos v,a\sinh u \sin v\big),$$ 
where $a$ is a positive number;
\vskip.04in 

\noindent {\rm (3)} A  flat surface  in $\mathbb E^4_2$  defined by 
\begin{align} \notag &\text{$ \frac{1}{\sqrt{2}m}$}\Big(\cosh(\sqrt{2}mx), \cosh(\sqrt{2}my), \sinh(\sqrt{2}mx),  \sinh(\sqrt{2}my)\Big),
\end{align} where $m$ is a positive number;
\vskip.04in 

\noindent {\rm (4)} A  flat surface  in $\mathbb E^4_2$  defined by 
$$\Big(0,\text{$ \frac{1}{a}$}\cosh (ax),\text{$ \frac{1}{a}$}\sinh(ax),y\Big),$$
where $a$ is a positive number;
\vskip.04in 

\noindent {\rm (5)} A  flat surface  in $\mathbb E^4_2$ defined by
\begin{align} \notag &\text{$ \(\frac{\cosh(\sqrt{2}x)}{\sqrt{2mr}},\frac{\cosh(\sqrt{2}y)}{\sqrt{2m(2m-r)}}, \frac{\sinh(\sqrt{2}x)}{\sqrt{2mr}}, \frac{\sinh(\sqrt{2}y)}{\sqrt{2m(2m-r)}} \)$},
\end{align} where $m$ and $r$ are positive numbers satisfying $2m>r>0$;

\noindent {\rm (6)} A  surface of negative curvature $-b^2$ in $\mathbb E^4_2$ defined by
\begin{equation}\begin{aligned} \notag & \Bigg(\text{$  \frac{1}{b}$}\cosh(bx)\cosh (by), \int_0^y\cosh (by)\sinh\!\(\! \text{$ \frac{4\sqrt{m^2-b^2}}{b}$} \tan^{-1}\!\Big( \tanh \text{$\frac{by}{2}$}\Big)\!\)dy, \\& \hskip.1in 
\text{$ \frac{1}{b}$}\sinh(bx)\cosh (by),\int_0^y\cosh (by)\cosh\!\( \!\text{$ \frac{4\sqrt{m^2-b^2}}{b}$} \tan^{-1}\!\Big( \tanh \text{$\frac{by}{2}$}\Big)\!\)dy\Bigg),\end{aligned}\end{equation}
where $b$ and $m$ are real numbers satisfying $0<b<m$.}

\section{Spacelike minimal surfaces with constant Gauss  curvature.} 

From the  equation of Gauss, we have

\vskip.1in
\noindent {\bf Lemma 7.1.}\label{L:7.1} {\it Let $M$ be a space-like minimal surface in $R^4_2(c)$. Then $K\geq c$. In particular, if $K=c$ holds identically, then $M$ is totally geodesic. }
\vskip.1in

For space-like {\it minimal surfaces} in  $R^4_2(c)$, Theorem 1 of \cite{S} implies that $M$ has constant Gauss curvature if and only if it has constant normal curvature.

We recall  the following result of M. Sasaki from \cite{S}.
  
\vskip.1in
\noindent {\bf Theorem 7.1.}\label{T:7.1} {\it Let $M$ be a space-like minimal surface in $R^4_2(c)$. If $M$ has constant Gauss curvature, then either 
  
 \vskip.04in
\noindent  {\rm (1)} $K=c$ and $M$ is a totally geodesic surface in $R^4_2(c)$;
  
\vskip.04in
\noindent   {\rm (2)} $c<0$, $K=0$  and $M$  is congruent to an open part of the  minimal surface defined by
  $$\text{$\frac{1}{\sqrt{2}}$}  \( \cosh u, \cosh v,0, \sinh u,  \sinh v \),$$ or   
  
\vskip.04in
\noindent   {\rm (3)} $c<0$, $K=c/3$ and $M$ is isotropic.}
  \vskip.1in

Let ${\bf R}^2$ be a plane with coordinates $s,t$.  Consider a map  ${\mathcal B}:{\bf R}^2\to \mathbb E^5_3$ given by
\begin{equation}\begin{aligned}&\label{7.1} {\mathcal B}(s,t)=\(\sinh \Big(\text{$\frac{2s}{\sqrt{3}}$}\Big)-\text{$ \frac{t^2}{3}$}-\(\text{$\frac{7}{8}+\frac{t^4}{18}$}\)e^{\frac{2s}{\sqrt{3}}}, t+\( \text{ $\frac{t^3}{3}-\frac{t}{4}$}\)e^{\frac{2s}{\sqrt{3}}},\right.
\\& \hskip1.1in 
 \text{ $\frac{1}{2}$} + \text{ $\frac{t^2}{2}$} e^{\frac{2s}{\sqrt{3}}},  t+ \(\text{ $\frac{t^3}{3}$}+ \text{ $\frac{t}{4}$}\)e^{\frac{2s}{\sqrt{3}}},\\&\hskip1.6in \left.
 \sinh \Big(\text{$\frac{2s}{\sqrt{3}}$}\Big)-\text{ $\frac{t^2}{3}$}-\text{$\(\frac{1}{8}+\frac{t^4}{18}\)$}e^{\frac{2s}{\sqrt{3}}}\). \end{aligned}\end{equation}
 The first author proved in \cite{c5} that ${\mathcal B}$ defines a full isometric parallel immersion  
  \begin{align}\label{7.2}\psi_{\mathcal B}:H^2(-\tfrac{1}{3})\to H^4_2(-1)\end{align}  
  of the hyperbolic plane $H^2(-\frac{1}{3})$ of curvature $-\frac{1}{3}$ into $H^4_2(-1)$. 
  
  The following result was also obtained in \cite{c5}.

\vskip.1in
\noindent {\bf Theorem 7.2.} \label{T:7.2} {\it Let $\psi:M\to H^4_2(-1)$ be a parallel full immersion of a space-like surface $M$ into $H^4_2(-1)$. Then    $M$ is minimal in $H^4_2(-1)$ if and only if  $M$ is congruent to an open part of the surface defined by
\begin{equation}\begin{aligned}&\notag \hskip.3in \(\sinh \Big(\text{$\frac{2s}{\sqrt{3}}$}\Big)-\text{$ \frac{t^2}{3}$}-\(\text{$\frac{7}{8}+\frac{t^4}{18}$}\)e^{\frac{2s}{\sqrt{3}}}, t+\( \text{ $\frac{t^3}{3}-\frac{t}{4}$}\)e^{\frac{2s}{\sqrt{3}}},\right.
\\& \hskip1in
 \text{ $\frac{1}{2}$} + \text{ $\frac{t^2}{2}$} e^{\frac{2s}{\sqrt{3}}},  t+ \(\text{ $\frac{t^3}{3}$}+ \text{ $\frac{t}{4}$}\)e^{\frac{2s}{\sqrt{3}}},\\& \hskip1.1in \left.
 \sinh \Big(\text{$\frac{2s}{\sqrt{3}}$}\Big)-\text{ $\frac{t^2}{3}$}-\text{$\(\frac{1}{8}+\frac{t^4}{18}\)$}e^{\frac{2s}{\sqrt{3}}}\). \end{aligned}\end{equation}}
 \vskip.1in

Combining Theorem 7.1 and Theorem 7.2, we obtain the following.

\vskip.1in
\noindent {\bf Theorem 7.3}\label{T:7.3} {\it Let $M$ be a non-totally geodesic space-like minimal surface in $H^4_2(-1)$. If $M$ has constant Gauss curvature $K$, then either 
  
 \vskip.04in
\noindent   {\rm (1)}  $K=0$ and $M$  is  congruent to an open part of the surface defined by
  $$\text{\small$\frac{1}{\sqrt{2}}$}  \( \cosh u, \cosh v,0, \sinh u,  \sinh v \),$$ or   

 \vskip.04in  
\noindent  {\rm (2)} $K=-\frac{1}{3}$ and $M$ is is congruent to an open part of the surface defined by
\begin{equation}\begin{aligned}&\notag \hskip.3in \(\sinh \Big(\text{$\frac{2s}{\sqrt{3}}$}\Big)-\text{$ \frac{t^2}{3}$}-\(\text{$\frac{7}{8}+\frac{t^4}{18}$}\)e^{\frac{2s}{\sqrt{3}}}, t+\( \text{ $\frac{t^3}{3}-\frac{t}{4}$}\)e^{\frac{2s}{\sqrt{3}}},\right.
\\& \hskip1in
 \text{ $\frac{1}{2}$} + \text{ $\frac{t^2}{2}$} e^{\frac{2s}{\sqrt{3}}},  t+ \(\text{ $\frac{t^3}{3}$}+ \text{ $\frac{t}{4}$}\)e^{\frac{2s}{\sqrt{3}}},\\& \hskip1.1in \left.
 \sinh \Big(\text{$\frac{2s}{\sqrt{3}}$}\Big)-\text{ $\frac{t^2}{3}$}-\text{$\(\frac{1}{8}+\frac{t^4}{18}\)$}e^{\frac{2s}{\sqrt{3}}}\). \end{aligned}\end{equation}}
 \vskip.1in

\section{Wintgen ideal surfaces satisfying $K^D=-2K$. }

We need the following existence result from \cite{CS}.

\vskip.1in
\noindent {\bf Theorem 8.1.}\label{T:8.1} {\it Let $c$ be a real number and  $\g$ with $3\g^2>-c$ be a positive solution of the second order partial differential equation
\begin{equation}\begin{aligned}\label{8.1}& \frac{\partial}{\partial x}\!\(\!\text{$\frac{(3\g \sqrt{c+3\g^2}-c)(6\g+2\sqrt{3c+9\g^2} \,)^{\sqrt{3}}\g_x}{2\g (c+3\g^2)}$}\)\\&
\hskip.2in  -\frac{\partial}{\partial y}\!\(\!\text{$\frac{(3\g \sqrt{c+3\g^2}-c)\g_y}{2\g (c+3\g^2)(6\g+2\sqrt{3c+9\g^2} \,)^{\sqrt{3}}}$}\)=\g\sqrt{c+3\g^2}\end{aligned}\end{equation}
 defined on a simply-connected domain $D\subset {\bf R}^2$.
 Then $M_{\g}=(D,g_{\g})$  with the metric 
\begin{align}\label{8.2}& g_\g=\text{$\frac{\sqrt{c+3\gamma^2} }{\gamma(6\gamma +2\sqrt{3c+9\gamma^2}\,)^{\sqrt{3}}}$}\!\(dx^2+(6\gamma +2\sqrt{3c+9\gamma^2}\,)^{2\sqrt{3}} dy^2\)\end{align}
 admits a non-minimal Wintgen ideal  immersion $\psi_\g: M_\g \to R^4_2(c)$  into a complete simply-connected indefinite space form $R^4_2(c)$  satisfying $K^D=2K$ identically. }
 \vskip.1in
 
The following result from \cite{CS} classifies  Wintgen ideal surfaces in $R^4_2(c)$ satisfying $K^D=2K$.

\vskip.1in
\noindent {\bf Theorem 8.2.}\label{T:8.2} {\it Let $M$ be a Wintgen ideal surface in a complete simply-connected indefinite space form $R^4_2(c)$ with $c=1,\,0$ or $-1$. If $M$ satisfies $K^D=2K$ identically, then  one of following three cases occurs:
  
 \vskip.04in
\noindent   {\rm (1)}  $c=0$ and $M$  is a totally geodesic surface in $\mathbb E^4_2$;
   
  \vskip.04in  
\noindent  {\rm (2)} $c=-1$  and $M$ is a minimal surface in $H^4_2(-1)$ congruent to an open part of $\psi_{\mathcal B}:H^2(-\tfrac{1}{3})\to H^4_2(-1)$ $\subset \mathbb E^5_3$  defined by
\begin{equation}\begin{aligned}&\notag \hskip.3in \(\sinh \Big(\text{$\frac{2s}{\sqrt{3}}$}\Big)-\text{$ \frac{t^2}{3}$}-\(\text{$\frac{7}{8}+\frac{t^4}{18}$}\)e^{\frac{2s}{\sqrt{3}}}, t+\( \text{ $\frac{t^3}{3}-\frac{t}{4}$}\)e^{\frac{2s}{\sqrt{3}}},\right.
\\& \hskip.1in \left.
 \text{ $\frac{1}{2}$} + \text{ $\frac{t^2}{2}$} e^{\frac{2s}{\sqrt{3}}},  t+ \(\text{ $\frac{t^3}{3}$}+ \text{ $\frac{t}{4}$}\)e^{\frac{2s}{\sqrt{3}}},
 \sinh \Big(\text{$\frac{2s}{\sqrt{3}}$}\Big)-\text{ $\frac{t^2}{3}$}-\text{$\(\frac{1}{8}+\frac{t^4}{18}\)$}e^{\frac{2s}{\sqrt{3}}}\); \end{aligned}\end{equation}
 
   \vskip.04in
\noindent   {\rm (3)} $M$ is a non-minimal surface in $R^4_2(c)$ which is congruent to an open part of $\psi_\g: M_\g\to R^4_2(c)$ associated with a positive solution $\g$ of the partial differential equation \e{8.1} as described in Theorem 8.1.}

\vskip.1in

\section{An application to minimal surfaces in $H^4_2(-1)$.}

A function $f$ on a space-like surface $M$ is called {\it   logarithm-harmonic}, if  $\Delta (\ln f)=0$ holds identically on $M$, where $\Delta (\ln f):=*d*(\ln f)$ is the Laplacian of $\ln f$ and $*$ is the Hodge star operator. A function $f$ on $M$ is called {\it subharmonic} if $\Delta f\geq 0$ holds everywhere on $M$.

In this section we present some results from \cite{c2010}.

\vskip.1in
\noindent {\bf Theorem 9.1.}\label{T:9.1} {\it Let $\psi:M\to H^4_2(-1)$ be a non-totally geodesic, minimal immersion of a space-like surface $M$ into $H^4_2(-1)$. Then  
\begin{align}\label{9.1} K+K^D\geq  -1\end{align} holds identically on $M$.

If $K+1$ is   logarithm-harmonic,  then the equality sign of \e{9.1} holds identically if and only if  $\psi:M\to H^4_2(-1)$ is congruent to an open portion of the immersion $\psi_{\phi}:H^2(-\frac{1}{3})\to H^4_2(-1)$ which is induced from the map $\phi:{\bf R}^2\to \mathbb E^5_3$ defined by
  \begin{equation}\begin{aligned}&\label{9.2} \phi (s,t)=\Bigg(\sinh \Big(\text{$\frac{2s}{\sqrt{3}}$}\Big)-\text{$ \frac{t^2}{3}$}-\(\text{$\frac{7}{8}+\frac{t^4}{18}$}\)e^{\frac{2s}{\sqrt{3}}}, t+\( \text{ $\frac{t^3}{3}-\frac{t}{4}$}\)e^{\frac{2s}{\sqrt{3}}},\\& \hskip1in 
 \text{ $\frac{1}{2}$} + \text{ $\frac{t^2}{2}$} e^{\frac{2s}{\sqrt{3}}},  t+ \(\text{ $\frac{t^3}{3}$}+ \text{ $\frac{t}{4}$}\)e^{\frac{2s}{\sqrt{3}}},\\&\hskip1.2in
 \sinh \Big(\text{$\frac{2s}{\sqrt{3}}$}\Big)-\frac{t^2}{3}-\text{$\(\frac{1}{8}+\frac{t^4}{18}\)$}e^{\frac{2s}{\sqrt{3}}}\Bigg). \end{aligned}\end{equation}}
 \vskip.1in

\noindent {\bf Corollary 9.1.}\label{C:9.1} {\it Let $\psi:M\to H^4_2(-1)$ be a minimal immersion of a space-like surface $M$ of constant Gauss curvature into $H^4_2(-1)$. 
Then the equality sign of \e{9.1} holds identically if and only if  one of the following two statements holds.

\vskip.04in
\noindent {\rm (1)} $K=-1, K^D=0,$ and $\psi$ is totally geodesic.

\vskip.04in
\noindent {\rm (2)}  $K^D=2K=-\frac{2}{3}$ and $\psi$ is congruent to an open part of the minimal surface $\psi_{\phi}:H^2(-\frac{1}{3})\to H^4_2(-1)$ induced from  \e{9.2}.}

\vskip.1in
\noindent {\bf Proposition 9.1.}\label{P:9.1} {\it Let $\psi:M\to {\mathbb E}^4_2$ be a minimal immersion of a space-like surface $M$  into the pseudo-Euclidean 4-space $\mathbb E^4_2$. Then  \begin{align}\label{9.3} K\geq -K^D\end{align} holds identically on $M$.

If $M$ has constant Gauss curvature, then 
 the equality sign of \e{9.3} holds identically if and only if $M$ is a totally geodesic surface.}

\vskip.1in
\noindent {\bf Proposition 9.2.}\label{P:9.2} {\it Let $\psi:M\to {\mathbb E}^4_2$ be a minimal immersion of a space-like surface $M$  into $\mathbb E^4_2$. We have

\vskip.04in
\noindent {\rm (1)}  If the equality sign of \e{9.1} holds identically, then $K$ is a non-logarithm-harmonic function.

\vskip.04in
\noindent {\rm (2)} If $M$ contains no totally geodesic points and  the equality sign of \e{9.3} holds identically on $M$, then $\ln K$ is  subharmonic.}
\vskip.1in
S
\vskip.1in
\noindent {\bf Proposition 9.3.}\label{P:9.3} {\it Let $\psi:M\to S^4_2(1)$ be a minimal immersion of a space-like surface $M$ into the neutral pseudo-sphere $S^4_2(1)$. Then  \begin{align}\label{9.4} K+K^D\geq 1\end{align} holds identically on $M$.

If $M$ has constant Gauss curvature, then the equality sign of \e{9.4} holds identically if and only if $M$ is a totally geodesic surface.}
\vskip.1in

Moreover, we have the following result from \cite{c2010}.

\vskip.1in
\noindent {\bf Proposition 9.4.}\label{P:9.4} {\it Let $\psi:M\to S^4_2(1)$ be a minimal immersion of a space-like surface $M$ into $S^4_2(1)$. We have

\vskip.04in
\noindent {\rm (1)} If the equality sign of \e{9.4} holds identically, then $K-1$ is non-logarithm-harmonic.

\vskip.04in
\noindent {\rm (2)} If $M$ contains no totally geodesic points and if the equality case of \e{9.4} holds, then $\ln (K-1)$ is subharmonic.}

\section{Wintgen ideal submanifolds are Chen submanifolds}

Consider a submanifold $M^n$ of a real space form $R^{n+m}(\epsilon)$, the normalized normal scalar curvature $\rho^\perp$ is defined as
$$\rho^\perp=\frac{2}{n(n-1)}\sqrt{\sum_{1\leq i<j\leq n; 1\leq r<s\leq m}\hskip-.3in \<R^\perp(e_i,e_j)\xi_r,\xi_s\>^2},$$
where $R^{\perp}$ is the normal connection of $M^{n}$, and $\{e_{1},\ldots, e_{n}\}$ and $\{\xi_{1},\ldots,\xi_{m}\}$ are the orthonormal frames of the tangent and normal bundles of $M^{n}$, respectively.

in 1999,  De Smet, Dillen, Vrancken and Verstraelen proved in \cite{DDVV} the
Wintgen inequality  \begin{align}\label{10.1} \rho\leq H^{2}-\rho^{\perp} +c\end{align} 
for all submanifolds $M^{n}$ of codimension 2 in all real space forms $\tilde M^{n+2}(c)$, where $\rho$ is the normalized scalar curvature defined by
\begin{align} \rho=\frac{2}{n(n-1)}\sum_{i<j}\<R(e_{i},e_{j})e_{j},e_{i}\>
\end{align} and $R$ is the Riemann curvature tensor of $M^{n}$.

The Wintgen inequality \e{10.1} was conjectured by  De Smet, Dillen, Vrancken and Verstraelen to hold for all submanifolds in all real space forms in the same paper \cite{DDVV}, known as {\it DDVV conjecture}.

Recently, Z. Lu \cite{Lu} and J. Ge and Z. Tang \cite{GT}, settled this conjecture independently in general. 
A submanifold $M^{n}$ of a real space form $M^{m}(c)$ is called a Wintgen ideal submanifold if it satisfies the equality case of \e{10.1} identically.

An $n$-dimensional submanifold $M$ of a  Riemannian manifold is called a {\it Chen submanifold} if 
\begin{align}\label{10.3}\sum_{i,j} \<\right. h(e_i,e_j),\overrightarrow{H}\left.\>h(e_i,e_j)\end{align}
 is parallel to the mean curvature vector
$\overrightarrow{H}$, where $h$ is the second fundamental form and $\{e_i\}$ is an orthonormal frame of the submanifold $M$ (cf. \cite{GVV}).

The following theorem was proved  by S. Decu, M. Petrovi\'c-Torga\v{s}ev and L. Verstraelen in \cite{dpv} which provides a very simple relationship between Wintgen ideal submanifolds and Chen submanifolds for submanifolds in real space forms.

\vskip.1in
\noindent  {\bf Theorem 10.1.} {\it Every Wintgen ideal submanifold of arbitrary dimension and codimension in a real space form is a Chen submanifold.}

\hfill{\vbox{\hbox{Department of Mathematics,}
             \hbox{Michigan State University,}
             \hbox{East Lansing, Michigan 48824,} \hbox{U.S.A.}
             \hbox{E-mail: {\tt bychen@math.msu.edu}}}

\end{document}